\documentclass[11pt]{article}
\usepackage[utf8]{inputenc}
\usepackage[english]{babel}
\usepackage[dvipsnames]{xcolor}
\usepackage{amsfonts,comment}
\usepackage{multicol}
\usepackage{amssymb,latexsym,amsmath}
\makeatletter \oddsidemargin 0in \evensidemargin 0in \textwidth
16cm \RequirePackage[dvips]{graphicx} \textheight 20cm
\DeclareGraphicsExtensions{.jpg, .eps}
\newtheorem{p1}{Proposition}[section]
\newtheorem{l1}{Lemma}[section]

\newtheorem{r1}{Remark}[section]

\definecolor{cadmiumred}{rgb}{0.89,0.0,0.13}

\begin{document}
\title{\bf
Errors due to departure from independence in multivariate Weibull distributions}
{\author{Subarna Bhattacharjee$^{1}$$^{}$\thanks{Corresponding author~:~~ Tel.: +91-9438004182,
E-mail: subarna.bhatt@gmail.com}, Aninda K. Nanda$^{2}$, Subhasree Patra$^{3}$ \\
{\it $^{1,3}$ Department of Mathematics,Ravenshaw University, Cuttack-753003, Odisha, India }\\
{\it $^{2}$ Department of Statistics, Babasahev BhimRao Ambedkar University,
Lucknow, UP, India}\\}
\date{}
\maketitle
\begin{abstract}
We do the error analysis in reliability measures due to the assumption of independence amongst the component lifetimes. In reliability theory, we come across different n-component structures like series, parallel, and k-out-of-n systems. A n component series system works only if all the n components work. While studying the reliability measures of a n-component series system, we mostly assume that all the components have independent lifetimes. Such an assumption eases mathematical complexity while analyzing the data and hence is very common. But in reality, the lifetimes of the components are very much interdependent. Such an assumption of independence hence leads to inaccurate analysis of data. In multiple situations like studying a complex system with many components, we turn to assuming independence keeping some room for error. However, if we have some knowledge of the behaviour of errors or some estimate on the error bound, we could decide if we assume independence and prefer mathematical simplicity (if we know the error is within our allowed limit), or keep the mathematical complexity and get accurate results without assuming independence. We aim to find the relative errors in the reliability measures for a n-component series system.
\end{abstract}
\section{Introduction}
Lu and Bhattacharyya I\\
As an application, Hougaard(1986) used the above Weibull model on a data on tumour appearance in 50 litters of female rats. Each litter contained one drug treated and two control rats. The data was studied to find estimate of hazards in marginal distribution, maximized likelihood function and standard errors. The true maximum likelihood estimate of log relative risk was found to be 0.944 for the dependence model and 0.04 less for the independence model. The former model had standard error 0.327 whereas the later had 0.01 less. The likelihood ratio test statistic was estimated to be 7.99 for the former and 0.29 more for the later model. It was also concluded that the two stage estimate and the maximum likelihood estimate differed marginally.\\
Crowder\\
Fitting the above model to a data on effect of lead absorption in small quantities into the bodies of rats for three lead level groups with dose levels 0, 0.1, 0.2, 0.4 and 0.8 each, Crowder(1989) studied the maximum log likelihoods for both linear and quadratic models, unconstraint and putting constraints. The measurements were noted as a pair of dry run and test run to analyse the drug effects for all the groups. The estimates varied between 0.40 and 0.71, 0.50 and 0.84, 0.71 and 0.84 for groups 1, 2 and 3 respectively with no constraint on the parameter l. If parameter $\gamma$ was taken unconstraint, the estimates varied between 0.04 and 1.46 for all the groups. For the unconstraint model, the log likelihood ratio test for linear and quadratic models gave $\chi_1^2$ estimate as 28.71, 7.03 and 2.77 for the three groups respectively. For the quadratic model the $\chi_2^2$ values for the three groups were found to be 13.86, 3.35 and 4.32 respectively. The plots between $P(T_1 > t_1)$ and $P(T_2 > t_2)$  showed no marked residual dependence between them. It was also observed from the plots that the curves moved to the right as the doses increased in all the groups. The drug had similar effect on all the groups qualitatively but with increasing lead level, the pain reducing effect decreased.\

In reliability theory, we come across different n component structure like series, parallel and k-out-of-n system. A n component series system works only if all the n components work. While studying reliability measures of a n component series system, we mostly assume that all the components have independent lifetimes. Such an assumption eases mathematical complexity while analysing the data and hence is very common. But in reality, the lifetimes of the components are very much interdependent. Such an assumption of independence hence leads to inaccurate analysis of data. To avoid such error in results, instead of independent distribution, many authors have preferred a number of distributions like, Marshall and Olkin’s multivariate exponential, Multivariate Gumbel’s type 1, Marshall and Olkin’s multivariate Weibull (1967a), multivariate Crowder (1989), multivariate Lee (1979), Lu and Bhattacharyya (1990) I, Farlie-Gumbel-Morgenstern Weibull, Lu and Bhattacharyya (1990) II, and many others. In multiple situations like studying complex system with many components, we turn to assuming independence keeping some room for the error. However, if we have some knowledge on the behaviour of errors or some estimate on the error bound, we could decide if we assume independence and prefer mathematical simplicity (if we know the error is within our allowed limit), or keep the mathematical complexity and get accurate results without assuming independence.

Moeschberger and Klein (1984) have studied the consequences of departures from independence when the component lifetimes in a series system are exponentially distributed. They have assumed the joint distribution to follow a Gumbel bivariate exponential model.

Klein and Moeschberger (1986) studied the magnitude of the errors under similar assumption about the component lifetimes to have independent exponential distributions when in fact the lifetimes follow the bivariate exponential distribution of Marshall and Olkin (series or parallel systems) or that of Freund (parallel systems).

Klein and Moeschberger (1987) investigated the consequences of departures from independence when the component lifetimes in a series system are exponentially distributed. They studied the departures by assuming the joint distribution to follow either one of the three Gumbel bivariate exponential models, the Downton bivariate exponential model, or the Oakes bivariate exponential model. They also concluded that the amount of error incurred depends upon the correlation between lifetimes and the relative mean life of the two components.

Gupta and Gupta (1990) studied relative error in reliability measures such as the reliability function, the failure rate and the mean residual life under the erroneous assumption of independence when in fact lifetimes follow a bivariate exponential model. The behaviour of these errors has been discussed to examine their structure as a function of time.

Recently, Nanda et al (2022) focussed on relative errors in the various reliability measures incurred due to assumption of independence of components in a two-component series system when components actually follow bivariate exponential distributions. They studied relative errors in four reliability measures viz. reliability function, failure rate function, mean residual life function and reversed hazard rate function.
\section{Mathematical Aspects}\label{intro}
${}$\hspace{0.8cm}The commonly used metrics for reliability analysis are survival function $\overline{F}(\cdot)=1-F(\cdot),$ failure rate function $r(\cdot)=f(\cdot)/\overline{F}(\cdot)$, reversed hazard rate function $\mu(\cdot)=f(\cdot)/F(\cdot)$ and aging intensity function $L(\cdot)$ where $f(\cdot)$  and $F(\cdot)$ are the density function and distribution functions respectively. These metrics describe the aging phenomena of a system represented by a nonnegative random variable $X$ representing the lifetime of any system or unit. The aging intensity function is mathematically expressed as
$$L(x)=\frac{-xf(x)}{\overline{F}(x)\ln \overline{F}(x)}=\frac{xr(x)}{\int_{0}^{x}r(u)du}, x>0.$$\\
A system with a well defined statistical distribution is said to belong to the non-parametric aging class of increasing(decreasing) failure rate, increasing(decreasing) failure rate average denoted by IFR(DFR) and IFRA (DFRA) according as $r(x)$ or $\frac{1}{x}\int_{0}^{x}r(x)dx$ is increasing(decreasing) in $x.$ Other aging classes having increasing(decreasing) aging intensity function called as increasing(decreasing) aging intensity classes exist  in literature. The corresponding aging classes  are respectively denoted by IAI (DAI).\\
\hspace*{1cm}Usually we deal with systems comprising of components (units) having some structure, namely series, parallel, $k$-out-of-$n.$ The components may work independently or have some dependence among them. The analysis of incurred error in computation of reliability metrics due to assumption of independence of components of a system when the components are actually dependent has been given importance by various authors.  As a simple case, in the present work we take up the case of a $n$ component series system.\\\hspace*{1cm} Recently, Nanda et al. (2022) focussed on relative errors in the various reliability measures incurred due to assumption of independence of components in a two component series system when components actually follow bivariate exponential distributions. In literature, we find works of Gupta and Gupta (1990), Klein and Moeschberger (1986, 1987), Moeschberger and Klein (1984) in computation and analysis of relative errors.\\ \hspace*{1cm}To start with, we recall some of the preliminaries and notations to be used in the present study.  The survival function of the lifetime $T(=\min_{1\leq i\leq k} X_i)$ of an $k$ component series system is written as\begin{eqnarray}
\overline{F}(x)
&=&P(X_1>x, X_2>x,\ldots, X_k>x)\nonumber\\
\label{system}
&=&\overline{F}(x,x,\ldots,x)
\end{eqnarray}
where $\overline{F}$ is the joint survival function of $X_1,X_2,\ldots,X_k.$\\
\hspace*{1cm}The relative errors in survival function $E^{\overline{F}}(x)$, failure rate function $E^{r}(x)$, reversed hazard rate function $E^{\mu}(x)$, aging intensity function $E^{L}(x)$ at time $x$ are respectively obtained as $E^{\overline{F}}(x)=(\overline{F}_{D}(x)-\overline{F}_{I}(x))/\overline{F}_{I}(x),$ $E^{r}(x)=(r_{D}(x)-r_{I}(x))/r_{I}(x),$ $E^{\mu}(x)=(\mu_{D}(x)-\mu_{I}(x))/\mu_{I}(x),$ $E^{L}(x)=(L_{D}(x)-L_{I}(x))/L_{I}(x).$ Here, $D$ and $I$ represent systems having dependent and independent components respectively.\\
\hspace*{1cm}We highlight some results in the form of lemma before this section ends. First, we give a lemma from Nanda et al. (2022).
\begin{l1}
\label{lemma}
For $x>0,$ $$g(x)=\frac{\gamma}{\beta}\Big(\frac{e^{\beta x}-1}{e^{\gamma x}-1}\Big)-1,$$
is increasing (resp. decreasing) in $x,$ provided $\beta>\gamma$ (resp. $\beta<\gamma$).
\end{l1}
\hspace*{1cm}We give a lemma from Bhattacharjee (2022). Sunoj and Rasin (2018) have proved the lemma using the concept of quantile function.
\begin{l1}\label{theoremss}
	A random variable $X$ is $IFRA$ (resp. $DFRA$) if and only if its $AI$ function satisfies $L(x) \geqslant (resp. \;\leqslant) 1$ for $x>0.$
\end{l1}
\section{Multivariate Exponential}
In this study, we first consider an $k$ component series system with dependent components having Marshall and Olkin's multivariate exponential distribution. To compute the associated errors as discussed in earlier section, we necessarily look into a series system of $n$ independent components having exponential distributions. Henceforth, by a series system, we mean a series system with $k$ components.
\subsection{Independent exponential}
\label{indenp}
The survival function of $k$ independent random variables $X_1,X_2,\ldots,X_k$ having exponential distribution is given by \begin{equation}\label{inde}\overline{F}(x_1,x_2,\ldots,x_k)=\exp\Big\{-\sum_{i=1}^{k}\lambda_i x_i\Big\}.\end{equation}
The survival function, failure rate, reversed failure rate and aging intensity functions of the series system formed out of $k$ independent components each having exponential distribution with survival function given in (\ref{inde}) are obtained using (\ref{system}) as
$$\overline{F}_{I}(x)=\exp\Big\{-x\sum_{i=1}^{k}\lambda_i\Big\}, r_{I}(x)=\sum_{i=1}^{k}\lambda_i, \mu_{I}(x)=\frac{\sum_{i=1}^{k}\lambda_i}{\exp\Big\{x\sum_{i=1}^{k}\lambda_i \Big\}-1}, L_{I}(x)=1.$$ Here, the notation $I$ stands for independent exponentially distributed random variables. Since $r_{I}(x)=\sum_{i=1}^{k}\lambda_i,$ we observe that the resultant series system is  exponentially distributed with parameter $\sum_{i=1}^{k}\lambda_i.$
\subsection{Marshall and Olkin's Multivariate Exponential distribution ($MOME$)}
The resulting joint distribution of lifetimes $X_1,\ldots,X_k$ of the components having $MOME$ is given by
\begin{eqnarray}
\overline{F}(x_1,x_2,\ldots,x_k)&=&\exp\Big\{-\sum_{i=1}^{k}\lambda_i x_i-\sum_{i_1<}\sum_{i_2}\lambda_{i_1,i_2}\max(x_{i_1},x_{i_2})\nonumber\\&&-\sum_{i_1<}\sum_{i_2<}\sum_{i_3}\lambda_{i_1,i_2,i_3}\max(x_{i_1},x_{i_2},x_{i_3})-\nonumber\\
\label{momed4}&&\ldots-\lambda_{i_1,i_2,\ldots,i_k}\max(x_{i_1},x_{i_2},\ldots,x_{i_k})
\Big\}
\end{eqnarray}
The survival function, failure rate, reversed failure rate and aging intensity functions of a series system formed out of $k$ components having $MOME$ distribution is obtained from (\ref{momed4}) and are given by
$\overline{F}_{D}(x)=\exp(-\lambda x), r_{D}(x)=\lambda, \mu_{D}(x)=\frac{\lambda }{\exp\big(\lambda x\big)-1}, L_{D}(x)=~1,$
where \begin{equation}
\label{lambmome}
\lambda=\Big(\sum_{i=1}^{k}\lambda_i+\sum_{i_1<}\sum_{i_2}\lambda_{i_1,i_2}+\sum_{i_1<}\sum_{i_2<}\sum_{i_3}\lambda_{i_1,i_2,i_3}+\ldots+\lambda_{i_1,i_2,\ldots,i_k}\nonumber
\Big).
\end{equation}
Here, $r_{D}(x)=\lambda,$ signifies that the resultant series system is  exponentially distributed with parameter $\lambda.$\\
The relative error incurred in survival function, failure rate, reversed failure rate and aging intensity functions by assuming independence of components when the components are actually dependent in formation of a series system are respectively given by
$$E_{1}^{\overline{F}}(x)=\exp\Big\{-x\Big(\lambda-\sum_{i=1}^{k}\lambda_i\Big)\Big\}-1, E_{1}^{r}(x)=\Big\{\Big(\lambda-\sum_{i=1}^{k}\lambda_i\Big)\Big\}\Big\{\sum_{i=1}^{k}\lambda_i\Big\}^{-1},$$ $$E_{1}^{\mu}(x)=\Big(\frac{\lambda}{\sum_{i=1}^{k}\lambda_i}\Big)\Big\{\frac{\exp(x\sum_{i=1}^{k}\lambda_i)-1}{\exp(\lambda x)-1}\Big\}-1, E_{1}^{L}(x)=0.$$
The following remark may be noted to prove the upcoming proposition \ref{RHR1}.
\begin{r1}
\label{pro}
Note that, for $g$ as defined in Lemma \ref{lemma}, it follows that $g(x)<0$ for $\gamma>\beta.$ This is because, if  $\gamma>\beta$ then $g(x)$  is decreasing in $x,$ giving $g(x)<g(0).$ Here, $g(0)=0,$ as
\begin{equation}
\lim_{x\rightarrow 0} \frac{\gamma}{\beta}\Big(\frac{e^{\beta x}-1}{e^{\gamma x}-1}\Big)=1\nonumber
\end{equation}
\end{r1}
The proposition throws some light on relative error in reversed hazard rate.
\begin{p1}
\label{RHR1}
So for $MOME$ distribution, the relative error in reversed hazard rate is negative for all $x\geq 0$ and decreasing in $x.$
\end{p1}
${}$\hspace{0.8cm} We note that relative error is a non-negative constant for failure rate and aging intensity function. However, for survival function, relative error is negative and decreasing in $x$. Also, relative error in reversed hazard rate is negative and decreasing in $x$.  We conclude that for $MOME$ distribution the assumption of components being independent leads to over-assessment of survival function and reversed hazard rate.  On the other hand, hazard rate and aging intensity functions give rise to under-assessment due to the assumption of independence. \\                                                                                                                                                                                                                                                                                   ${}$\hspace{0.8cm} Moreover, the absolute error in survival function and reversed hazard rate are bounded between 0 and 1. The absolute error is found to be least in aging intensity function as expected as the function gives an overall idea of aging phenomenon of a system and is observed to be least affected by the assumption of independence of components, when the components are actually dependent.
\subsection{Multivariate Gumbel's Type 1 ($MG1$)}
The resulting joint distribution of lifetimes $X_1,\ldots,X_k$ of the components having $MG1$ is given by
\begin{eqnarray}
\overline{F}(x_1,x_2,\ldots,x_k)&=&\exp\Big\{-\sum_{i=1}^{k}\lambda_i x_i-\sum_{i_1<}\sum_{i_2}\lambda_{i_1,i_2}x_{i_1}x_{i_2}-\sum_{i_1<}\sum_{i_2<}\sum_{i_3}\lambda_{i_1,i_2,i_3}x_{i_1}x_{i_2}x_{i_3}-\ldots-\nonumber\\
&&\lambda_{i_1,i_2,\ldots,i_k}x_{i_1}x_{i_2}\ldots x_{i_k})
\Big\}\nonumber.
\end{eqnarray}
For simplification purpose we use the following notations $$a_1=\sum_{i=1}^{k}\lambda_i;~ a_2=\sum_{i_1<}\sum_{i_2}\lambda_{i_1,i_2};~ a_3=\sum_{i_1<}\sum_{i_2<}\sum_{i_3}\lambda_{i_1,i_2,i_3};~ \ldots,~ a_k=\lambda_{i_1,i_2,\ldots,i_k}.$$ The survival function, failure rate, reversed failure rate and aging intensity functions of a series system formed out of $k$ components having Gumbel I model are respectively given by
\begin{equation}
\overline{F}_{D}(x)=\exp\big\{-a_1x -a_2x^2-a_3x^3-\ldots-a_kx^k
\big\}\nonumber
\end{equation}
\begin{equation}
r_{D}(x)=a_1 +2a_2x+3a_3x^2+\ldots+ka_kx^{k-1}\nonumber
\end{equation}
\begin{eqnarray}
\mu_{D}(x)&=&\frac{a_1 +2a_2x+3a_3x^2+\ldots+ka_kx^{k-1}}{\exp\big\{a_1x+a_2x^2+a_3x^3+\ldots+a_kx^k
\big\}-1}\nonumber\\
&=&\frac{\theta^{'}(x)}{e^{\theta(x)}-1},\nonumber
\end{eqnarray}
where $$\theta(x)
=a_1x+a_{2}x^{2}+a_3x^{3}+\ldots+a_kx^{k},$$ and
$\theta^{'}(x)$ represents differentiation of $\theta(x)$ with respect to $x.$
\begin{equation}
L_{D}(x)=\frac{a_1x+2a_2 x^2+3a_3 x^3+\ldots+ka_k x^{k}
}{a_1x+a_2x^2+a_3x^3+\ldots+a_kx^k
}\nonumber
\end{equation}
From Lemma \ref{bhatt}, it follows that $1\leq L_{D} (x)\leq k$ and we conclude that the resultant series system is $IFRA.$\\
The relative errors so incurred in $MG1$ distributions are given by
$$E_{2}^{\overline{F}}(x)=\exp\big\{ -a_2x^2-a_3x^3-\ldots-a_lx^k
\big\}-1$$
$$E_{2}^{r}(x)=\frac{2a_2x+3a_3x^2+\ldots+ka_kx^{k-1}
}{a_0}$$
\begin{equation}
E_{2}^{\mu}(x)= \Big(\frac{a_0+\theta^{'}(x)}{a_0}\Big)\Big(\frac{e^{a_0x}-1}{e^{a_0x+\theta(x)}-1}\Big)-1\nonumber
\end{equation}
$$E_{2}^{L}(x)=\frac{a_1 x^2+2a_2x^3+\ldots+(k-1)a_k x^{k-1}
}{a_0x+a_1 x^2+a_2x^3+\ldots+a_k x^k}
$$
${}$\hspace{0.8cm} In Gumbel I model, relative error is negative  and decreasing in $x$ for survival function. The relative error in hazard rate is non-negative and increasing in $x$.  The relative error in reversed hazard rate. However, for survival function, relative error is  decreasing in $x$. We conclude that for $MG1$ distribution the assumption of components being independent leads to over-assessment of survival function and under-assessment. Using Lemma \ref{bhatt}, we find that an upper bound of relative error in $AI$ function  is $(k-1).$\\
\hspace*{1cm} In this section, we analyze the errors in computation of various reliability measures occurring due to the assumption of independent components instead of working with $k$ number of dependent components in a series system. We focussed on Marshall Olkin and Gumbel' s Type 1 multivariate exponential distributions. Similar study can be extended for other multivariate exponential distributions. Now we direct our study to multivariate Weibull distributions which are generalizations of multivariate exponential distribution.
\section{Multivariate  Weibull Models}
We now take up various multivariate Weibull distributions and do similar study as done in earlier section.  To start with, we essentially have to figure out the reliability measures for a series system comprising of $k$ independent
components having Weibull distribution.
\subsection{Independent Weibull}
The survival function of $k$ independent random variables $X_1,X_2,\ldots,X_k$ having  Weibull distribution is given by
\begin{equation}
\label{indwe}
\overline{F}(x_1,x_2,\ldots,x_k)=\exp\Big\{-\sum_{i=1}^{k}\lambda_i x_i^{\alpha_i}\Big\}.
\end{equation}
The survival function, the failure rate, the reversed failure rate and the aging intensity functions of a series system formed out of $n$ independent components each having Weibull distribution with survival function given in (\ref{indwe}) are obtained using (\ref{system}) as
$$\overline{F}_{I}(x)=\exp\Big\{-\sum_{i=1}^{k}\lambda_i x^{\alpha_i}\Big\}, r_{I}(x)=\sum_{i=1}^{k}\lambda_i \alpha_i x^{\alpha_i-1},$$ $$\mu_{I}(x)=\frac{\sum_{i=1}^{k}\lambda_i \alpha_i x^{\alpha_i-1}}{\exp\Big\{\sum_{i=1}^{k}\lambda_i x^{\alpha_i}\Big\}-1}, L_{I}(x)=\frac{\sum_{i=1}^{k}\lambda_i \alpha_i x^{\alpha_i}}{\sum_{i=1}^{k}\lambda_i x^{\alpha_i}}.$$ Henceforth, the notation $I$ is used for independent Weibull random variables. \\
\hspace*{1cm}Applying Lemma \ref{bhatt}, we find that $$\min_{i} \alpha_i\leq L_{I}(x)\leq \max_{i} \alpha_i$$ claim that the resultant system is $IFRA$ ($DFRA$) according as $\min \alpha_i>1 (\max\alpha_i<1).$
\subsection{Marshal and Olkin (1967a)(MOMW)}
Marshal and Olkin (1967a) and Lee and Thompson (1974) gave a generalization of multivariate Weibull distribution, called as Marshal and Olkin multivariate Weibull Models with survival function given by
\begin{eqnarray}
\overline{F}(x_1,x_2,\ldots,x_k)&=&\exp\Big\{-\sum_{i=1}^{k}\lambda_i x_i^{\alpha_i}-\sum_{i_1<}\sum_{i_2}\lambda_{i_1,i_2}\max(x_{i_1}^{\alpha_{i_1}},x_{i_2}^{\alpha_{i_2}})\nonumber\\
\label{momed}
&&-\sum_{i_1<}\sum_{i_2<}\sum_{i_3}\lambda_{i_1,i_2,i_3}\max(x_{i_1}^{\alpha_{i_1}},x_{i_2}^{\alpha_{i_2}},x_{i_3}^{\alpha_{i_3}})-\nonumber\\&&\ldots-\lambda_{i_1,i_2,\ldots,i_k}\max(x_{i_1}^{\alpha_{i_1}},x_{i_2}^{\alpha_{i_2}},\ldots,x_{i_k}^{\alpha_{i_k}})
\Big\}
\end{eqnarray}
The survival function, the failure rate, the reversed failure rate and the aging intensity functions of a series system formed out of $n$ components having $MOMW$ are respectively given by
$\overline{F}_{D}(x)=\exp(A(x)),$ $r_{D}(x)=A^{'}(x),$
$\mu_{D}(x)=\frac{A^{'}(x)}{\exp[A(x)]-1},$ $L_{D}(x)=\frac{xA^{'}(x)}{A(x)}$
where 
$$A(x)=\Big\{\sum_{i=1}^{k}\lambda_i x^{\alpha_i}+\sum_{i_1<}\sum_{i_2}\lambda_{i_1,i_2}~x^{\max(\alpha_{i_1}, \alpha_{i_2})}+\ldots+\lambda_{i_1,i_2,\ldots,i_k}x^{\max(\alpha_{i_1}, \alpha_{i_2},\ldots,\alpha_{i_k})}
\Big\}.$$
 The relative error incurred in survival function, failure rate, reversed failure rate and aging intensity functions by assuming independence of components where the components are actually dependent in formation of a series system are respectively given by
$$E_{3}^{\overline{F}}(x)=\exp\Big(-A(x)+\sum_{i=1}^{k}\lambda_i x^{\alpha_i}\Big)-1$$
$$E_{3}^{r}(x)=\frac{1}{x}\Big(A(x)-\sum_{i=1}^{k}\lambda_i x^{\alpha_i}\Big)\Big\{\sum_{i=1}^{k}\lambda_i \alpha_i x^{\alpha_i-1}\Big\}^{-1}$$
\subsection{Multivariate Crowder (1989) ($MC$)}
Hougaard (1986,1989) presented a multivariate  Weibull distribution with joint survival function
\begin{equation}
\overline{F}(x_1,x_2,\ldots,x_k)=\exp\Big\{-\big(\sum_{i=1}^{k}\lambda_i x_{i}^{\alpha_i}\big)^{l}\Big\},
\end{equation}
where $l>0, p \geq 0$ and $x_i \geq 0.$ Crowder (1989) extended Hougaard's distributions and proposed "multivariate distributions with Weibull connections" with
\begin{equation}
\overline{F}(x_1,x_2,\ldots,x_k)=\exp\Big\{\gamma^{l}-\big(\gamma+\sum_{i=1}^{k}\lambda_i x_{i}^{\alpha_i}\big)^{l}\Big\}
\end{equation}
where $l>0, \nu \geq 0$ and $\alpha_i >0.$ In the special case when $\alpha_1=\ldots=\alpha_k=\alpha,$ the marginals are all Weibull with the same parameter.\\
The survival function, failure rate, reversed failure rate and aging intensity functions of a series system formed out of $n$ components having Crowder distribution are respectively given by
\begin{equation}
\overline{F}(x)=\exp\Big\{\gamma^{l}-\big(\gamma+\sum_{i=1}^{k}\lambda_i x^{\alpha_i}\big)^{l}\Big\}\nonumber
\end{equation}
\begin{equation}
r(x)= l\Big(\gamma+\sum_{i=1}^{k}\lambda_i x^{\alpha_i}\Big)^{l-1}\Big(\sum_{i=1}^{k}\lambda_i \alpha_i x^{\alpha_i-1}\Big)\nonumber
\end{equation}
 $$\mu(x)=\frac{l\Big(\gamma+\sum_{i=1}^{k}\lambda_i x^{\alpha_i}\Big)^{l-1}\Big(\sum_{i=1}^{k}\lambda_i \alpha_i x^{\alpha_i-1}\Big)}{\exp\Big\{-\gamma^{l}+\big(\gamma+\sum_{i=1}^{k}\lambda_i x^{\alpha_i}\big)^{l}\Big\}-1},$$
  $$L(x)=\frac{l\Big(\gamma+\sum_{i=1}^{k}\lambda_i x^{\alpha_i}\Big)^{l-1}\Big(\sum_{i=1}^{k}\lambda_i \alpha_i x^{\alpha_i}\Big)}{\exp\Big\{\gamma^{l}-\big(\gamma+\sum_{i=1}^{k}\lambda_i x^{\alpha_i}\big)^{l}\Big\} \Big\{\gamma^{l}-\big(\gamma+\sum_{i=1}^{k}\lambda_i x^{\alpha_i}\big)^{l}\Big\}}$$
  The relative error incurred in survival function, failure rate, reversed failure rate and aging intensity functions by assuming independence of components where the components are actually dependent in formation of a series system are respectively given by
 $$E_{4}^{\overline{F}}(x)=\exp\big\{\gamma^{l}+\sum_{i=1}^{k} \lambda_{i}x^{\alpha_i}-\big(\gamma+\sum_{i=1}^{k} \lambda_{i}x^{\alpha_i}\big)^{l}\big\}-1,E_{4}^{r}(x)=l\Big(\gamma+\sum_{i=1}^{k} \lambda_{i}x^{\alpha_i}\Big)^{l-1}-1,$$
 $$E_{4}^{\mu}(x)=\frac{l\Big(\gamma+\sum_{i=1}^{k}\lambda_i x^{\alpha_i}\Big)^{l-1}}{\Big\{\exp\big\{-\gamma^{l}+\big(\gamma+\sum_{i=1}^{k}\lambda_i x^{\alpha_i}\big)^{l}\big\}-1\Big\}}\Big\{\exp\Big(\sum_{i=1}^{k}\lambda_i x^{\alpha_i}\Big)-1\Big\}.$$
 $$E_{4}^{L}(x)=\frac{l\Big(\gamma+\sum_{i=1}^{k}\lambda_i x^{\alpha_i}\Big)^{l-1}\Big(\sum_{i=1}^{k}\lambda_i x^{\alpha_i}\Big)}{\Big\{\exp\Big(\gamma^{l}-\big(\gamma+\sum_{i=1}^{k}\lambda_i x^{\alpha_i}\big)^{l}\Big)\Big\} \Big\{\gamma^{l}-\big(\gamma+\sum_{i=1}^{k}\lambda_i x^{\alpha_i}\big)^{l}\Big\}}-1.$$
In particular, Lee (1979): II ($MLII$)
\begin{equation}
\overline{F}(x_1,x_2,\ldots,x_k)=\exp\big\{-\big(\sum_{i=1}^{k} \lambda_{i}x_i^{\alpha_i}\big)^{l}\big\}
\end{equation}
where $\alpha_i>0, 0<\gamma\leq 1, \lambda_i>0,x_i\geq 0.$\\
\subsection{Multivariate Lee (1979) ($ML$)}
\begin{eqnarray}
\overline{F}(x_1,x_2,\ldots,x_k)&=&\exp\Big\{-\sum_{i=1}^{k}\lambda_i c_i^{\alpha} x_i^{\alpha}-\sum_{i_1<}\sum_{i_2}\lambda_{i_1,i_2}\max(c_{i_1}^{\alpha}x_{i_1}^{\alpha},c_{i_2}^{\alpha}x_{i_2}^{\alpha})\nonumber\\
\label{momed}
&&-\sum_{i_1<}\sum_{i_2<}\sum_{i_3}\lambda_{i_1,i_2,i_3}\max(c_{i_1}^{\alpha}x_{i_1}^{\alpha},c_{i_2}^{\alpha}x_{i_2}^{\alpha},c_{i_3}^{\alpha}x_{i_3}^{\alpha})-\nonumber\\
&&\ldots-\lambda_{i_1,i_2,\ldots,i_k}\max(c_{i_1}^{\alpha}x_{i_1}^{\alpha},c_{i_2}^{\alpha}x_{i_2}^{\alpha},\ldots,c_{i_k}^{\alpha}x_{i_k}^{\alpha})
\Big\}
\end{eqnarray}
Taking its independent counterpart as $\overline{F}(x_1,x_2,\ldots,x_k)=\exp\Big\{-\sum_{i=1}^{k}\lambda_i c_i^{\alpha} x_i^{\alpha}\Big\},$ we note that the survival function, failure rate function, reversed hazard rate and ageing intensity function of a series system with $k$ independent components are given by $\overline{F}_{I}(x)=\exp\Big\{-x^{\alpha}\sum_{i=1}^{k}\lambda_i c_i^{\alpha} \Big\},$ $r_{I}(x)=\alpha x^{\alpha-1}\sum_{i=1}^{k}\lambda_i c_{i}^{\alpha}, \mu_{I}(x)=\frac{\alpha x^{\alpha-1}\sum_{i=1}^{k}\lambda_i c_{i}^{\alpha}}{\exp\Big\{x^{\alpha}\sum_{i=1}^{k}\lambda_i c_{i}^{\alpha} \Big\}-1},$ and $L_{I}(x)=\alpha.$\\
\hspace*{1cm}The survival function, failure rate, reversed failure rate and aging intensity functions of a series system formed out of $n$ components having multivariate $ML$ distribution are respectively given by
$$\overline{F}(x)=\exp\Big\{-\lambda_{L}x^{\alpha}
\Big\}, r(x)=\alpha x^{\alpha-1}\lambda_{L},
\mu(x)=\frac{\alpha x^{\alpha-1}\lambda_{L}}{e^{\lambda_{L}x^{\alpha}}-1},
  L(x)=\alpha,$$ where \begin{eqnarray}
\lambda_{L}&=&\Big(\sum_{i=1}^{k}\lambda_i c_i^{\alpha}+\sum_{i_1<}\sum_{i_2}\lambda_{i_1,i_2}\max(c_{i_1}^{\alpha},c_{i_2}^{\alpha})+\sum_{i_1<}\sum_{i_2<}\sum_{i_3}\lambda_{i_1,i_2,i_3}\max(c_{i_1}^{\alpha},c_{i_2}^{\alpha},c_{i_3}^{\alpha})+\nonumber\\
&&\ldots+\lambda_{i_1,i_2,\ldots,i_k}\max(c_{i_1}^{\alpha},c_{i_2}^{\alpha},\ldots,c_{i_k}^{\alpha})\Big).
\end{eqnarray}
A little modification of Lemma \ref{lemma} gives the following result.
\begin{l1}
\label{lemma2}
For $x>0,$ $\beta,\gamma,\alpha>0,$ $$h(x)=\frac{\gamma}{\beta}\Big(\frac{e^{\beta x^{\alpha}}-1}{e^{\gamma x^{\alpha}}-1}\Big)-1,$$
is increasing (resp. decreasing) in $x,$ provided $\beta>\gamma$ (resp. $\beta<\gamma$).
\end{l1}
{\bf Proof.} Taking $x^{\alpha}=t (say),$ and noting that $$\frac{d}{dx}h(x)=\frac{dh}{dt}\frac{dt}{dx},$$ we complete the proof.\\
\begin{r1}
\label{pro}
Note that, for $h$ as defined in Lemma \ref{lemma2}, $h(x)<0$ for $\gamma>\beta.$ This is because, if  $\gamma>\beta$ then $h(x)$  is decreasing in $x,$ giving $h(x)<h(0).$ Here, $h(0)=0,$ as
\begin{equation}
\lim_{x\rightarrow 0} \frac{\gamma}{\beta}\Big(\frac{e^{\beta x^{\alpha}}-1}{e^{\gamma x^{\alpha}}-1}\Big)=1\nonumber
\end{equation}
\end{r1}
 The relative error incurred in survival function, failure rate, reversed failure rate and aging intensity functions by assuming independence of components where the components are actually dependent in formation of a series system are respectively given by
$$E_{6}^{\overline{F}}(x)=\exp\Big\{-x^{\alpha}\Big(\lambda_{L}-\sum_{i=1}^{k}\lambda_{i}c_{i}^{\alpha}\Big)\Big\}-1,
E_{6}^{r}(x)=\frac{\lambda_{L}}{\sum_{i=1}^{k}\lambda_{i}c_{i}^{\alpha}}-1,$$ $$E_{6}^{\mu}(x)=\Big(\frac{\lambda_{L}}{\sum_{i=1}^{k}\lambda_ic_{i}^{\alpha}}\Big)\Big\{\frac{\exp(x^{\alpha}\sum_{i=1}^{k}\lambda_i c_{i}^{\alpha})-1}{\exp(\lambda_{L} x^{\alpha})-1}\Big\}-1,E^{L}_{6}(x)=0.$$
We observe that the relative error in  each of $sf,fr,rhr$ is negative, whereas in $ai$ it is zero.
\subsection{Lu and Bhattcharyya (1990): I}
\begin{equation}
\overline{F}(x_1,x_2,\ldots,x_k)=\exp\Big\{-\Big(\sum_{i=1}^{k}\lambda_{i}x_i^{\alpha_i}+\delta w(x_1,x_2,\ldots,x_k)\Big)\Big\}
\end{equation}
where $$w(x_1,x_2,\ldots,x_k)=\Big\{\sum_{i=1}^{k}\lambda_{i}^{\frac{1}{m}}x_i^{\frac{\alpha_i}{m}}\Big\}^{m}$$
The survival function, failure rate, reversed failure rate and aging intensity functions of a series system formed out of $k$ components having Crowder distribution are respectively given by
\begin{equation}
\overline{F}^{T}_{1}(x)=\exp\Big\{-\Big\{\sum_{i=1}^{k}\lambda_{i}x^{\alpha_i}+\delta \Big(\sum_{i=1}^{k}\lambda_{i}^{\frac{1}{m}}x^{\frac{\alpha_i}{m}}\Big)^{m}\Big\}\Big\}\nonumber
\end{equation}
\begin{equation}
r^{T}_1(x)= \sum_{i=1}^{k}\lambda_{i}\alpha_i x_i^{\alpha_i}+\frac{\delta}{x}\Big(\sum_{i=1}^{k}\lambda_{i}^{\frac{1}{m}}\alpha_i x^{\frac{\alpha_i}{m}}\Big)\Big(\sum_{i=1}^{k}\lambda_{i}^{\frac{1}{m}}x^{\frac{\alpha_i}{m}}\Big)^{m}\nonumber
\end{equation}

 $$\mu^{T}_1(x)=\frac{\gamma \big(\sum_{i=1}^{k} \lambda_{i}x^{\alpha_i}\big)^{\gamma-1} \Big(\sum_{i=1}^{k}\lambda_i \alpha_i x^{\alpha_i-1}\Big)}{\exp\Big\{\sum_{i=1}^{k}\lambda_i x^{\alpha_i}\Big\}-1},$$

  $$L^{T}_1(x)=\frac{\Big(\sum_{i=1}^{k}\lambda_{i}\alpha_i x_i^{\alpha_i}\Big)+\delta \Big(\sum_{i=1}^{k}\lambda_{i}^{\frac{1}{m}}\alpha_i x^{\frac{\alpha_i}{m}}\Big)\Big(\sum_{i=1}^{k}\lambda_{i}^{\frac{1}{m}} x^{\frac{\alpha_i}{m}}\Big)^{m-1}}{\Big(\sum_{i=1}^{k}\lambda_{i} x_i^{\alpha_i}\Big)+\delta \Big(\sum_{i=1}^{k}\lambda_{i}^{\frac{1}{m}} x^{\frac{\alpha_i}{m}}\Big)^{m}}$$
%
%
%
\section{Analysis of parallel system with $n$ components}
\begin{equation*}
	\begin{array}{rll}
		F(x_1,x_2,x_3,\cdots,x_n)=&P(X_1\leq x_1,X_2\leq x_2,X_3\leq x_3,\cdots, X_n\leq x_n)&\\
		=&P(\cap_{i=1}^n X_i\leq x_i)&\\
		=&1-P(\cup_{i=1}^n X_i>x_i)&\text{(By De Morgan's Law)}\\
		=&1-\bigg[\sum_{i=1}^n P(X_i> x_i)&\\
		&-\underset{1\leq i<j\leq n}{\sum\sum}P(X_i> x_i, X_j> x_j)&\\
		&+\underset{1\leq i<j<k\leq n}{\sum\sum\sum}P(X_i> x_i, X_j> x_j, X_k> x_k)&\\
		&\vdots&\\
		&(-1)^{n-1}P(X_1> x_1,X_2> x_2,\cdots,X_n> x_n)\bigg]&\text{(By Poincare's Theorem)}\\
		=&1-\sum_{i=1}^n \bar{F}_{X_i}(x_i)+\underset{1\leq i<j\leq n}{\sum\sum}\bar{F}_{X_i,X_j}(x_i,x_j)&\\
		&-\underset{1\leq i<j<k\leq n}{\sum\sum\sum}\bar{F}_{X_i,X_j,X_k}(x_i,x_j,x_k)\\
		&\vdots&\\
		&(-1)^n\bar{F}(x_1,x_2,x_3,\cdots,x_n)&\\
		=&1-\sum_{i=1}^n \bar{F}(-\infty,\cdots,-\infty,x_i,-\infty,\cdots,-\infty)&\\
		&\multicolumn{2}{l}{+\underset{1\leq i<j\leq n}{\sum\sum}\bar{F}(-\infty,\cdots,-\infty,x_i,-\infty,\cdots,-\infty,x_j,-\infty,\cdots,-\infty)}\\
		&\multicolumn{2}{l}{-\underset{1\leq i<j<k\leq n}{\sum\sum\sum}\bar{F}(-\infty,\cdots,-\infty,x_i,-\infty,\cdots, -\infty,x_j,-\infty,\cdots,-\infty,x_k,-\infty,\cdots,-\infty)}\\
		&\vdots&\\
		&(-1)^n\bar{F}(x_1,x_2,x_3,\cdots,x_n)&\\
	\end{array}
\end{equation*}
\begin{equation*}
	\begin{array}{rll}
		\bar{F}_{T}(t)=1-F(t,t,t,\cdots,t)=&\sum_{i=1}^n \bar{F}(0,\cdots,0,t,0,\cdots,0)&\\
		&\multicolumn{2}{l}{-\underset{1\leq i<j\leq n}{\sum\sum}\bar{F}(0,\cdots,0,t,0,\cdots,0,t,0,\cdots,0)}\\
		&\multicolumn{2}{l}{+\underset{1\leq i<j<k\leq n}{\sum\sum\sum}\bar{F}(0,\cdots,0,t,0,\cdots, 0,t,0,\cdots,0,t,0,\cdots,0)}\\
		&\vdots&\\
		&(-1)^{n-1}\bar{F}(t,t,t,\cdots,t)&\\
	\end{array}
\end{equation*}

\subsection{Independent exponential}
\begin{equation}\label{inde}\overline{F}(x_1,x_2,\ldots,x_n)=\exp\Big\{-\sum_{i=1}^{n}\lambda_i x_i\Big\}.
\end{equation}
\begin{equation*}
	\begin{split}
		\bar{F}_{T}(t)=&\sum_{i=1}^n\exp\left\{-t\lambda_i\right\}\\
		&-\mathop{\sum\sum}_{1\leq i<j\leq n}\exp\left\{-t(\lambda_i+\lambda_j)\right\}\\
		&+\mathop{\sum\sum\sum}_{1\leq i<j<k\leq n}\exp\left\{-t(\lambda_i+\lambda_j+\lambda_k)\right\}\\
		&\vdots\\
		&+(-1)^{n-1}\exp\left\{-t\sum_{i=1}^{n}\lambda_i\right\}\\
		=&\sum_{k=1}^n(-1)^{k-1}\mathop{\sum\cdots\sum}_{1\leq i_1<\cdots<i_k\leq n}\exp\left\{-t\left(\sum_{p=1}^k\lambda_{i_p}\right)\right\}
	\end{split}
\end{equation*}

\subsection{Marshall and Olkin's Multivariate Exponential distribution ($MOME$)}
The resulting joint distribution of lifetimes $X_1,\ldots,X_n$ of the components having $MOME$ is given by
\begin{eqnarray}
\overline{F}(x_1,x_2,\ldots,x_n)&=&\exp\Big\{-\sum_{i=1}^{k}\lambda_i x_i-\sum_{i_1<}\sum_{i_2}\lambda_{i_1,i_2}\max(x_{i_1},x_{i_2})\nonumber\\&&-\sum_{i_1<}\sum_{i_2<}\sum_{i_3}\lambda_{i_1,i_2,i_3}\max(x_{i_1},x_{i_2},x_{i_3})-\nonumber\\
\label{momed4}&&\ldots-\lambda_{1,2,\ldots,n}\max(x_{1},x_{2},\ldots,x_{n})
\Big\}
\end{eqnarray}
\begin{equation*}
	\begin{split}
		\bar{F}_{T}(t)=&\sum_{i=1}^n\exp\left\{-t\lambda_i-\sum_{j\neq i}t\lambda_{i,j}-\mathop{\mathop{\sum\sum}_{j\neq i\neq k}}_{j<k}t\lambda_{i,j,k}-\cdots-t\lambda_{1,\cdots,k}\right\}\\
		-&\mathop{\sum\sum}_{1\leq i_1<i_2\leq n}\exp\left\{-t\lambda_{i_1}-t\lambda_{i_2}-t\lambda_{i_1,i_2}-\sum_{j\notin\{i_1,i_2\}}\left(t\lambda_{i_1,j}+t\lambda_{i_2,j}\right)-\sum_{j\notin\{i_1,i_2\}}t\lambda_{i_1,i_2,j}\right.\\
		&\left. -\mathop{\sum\sum}_{j<k\notin\{i_1,i_2\}}\left(t\lambda_{i_1,j,k}+t\lambda_{i_2,j,k}\right)-\cdots-t\lambda_{1,2,\cdots,n}\right\}\\
		\vdots&\\
		+&(-1)^{n-1}\exp\left\{-t\left(\sum_{i=1}^n\lambda_i +\mathop{\sum\sum}_{i<j}\lambda_{i,j}+\mathop{\sum\sum\sum}_{i<j<k}\lambda_{i,j,k}+\cdots+\lambda_{1,2,\cdots,n}\right)\right\}\\
	\end{split}
\end{equation*}
\begin{equation*}
	\begin{split}
		\bar{F}_T(t)=&\sum_{i=1}^n\exp\left\{-t\left(\lambda_i+\sum_{j\neq i}\lambda_{i,j}+\mathop{\mathop{\sum\sum}_{j\neq i\neq k}}_{j<k}\lambda_{i,j,k}+\cdots+\lambda_{1,2,\cdots,n}\right)\right\}\\
		-&\mathop{\sum\sum}_{1\leq i_1<i_2\leq n}\exp\left\{-t\left(\lambda_{i_1}+\lambda_{i_2}+\lambda_{i_1,i_2}+\sum_{j\notin\{i_1,i_2\}}(\lambda_{i_1,j}+\lambda_{i_2,j}+\lambda_{i_1,i_2,j})\phantom{\mathop{\mathop{\sum\sum}_{j,k\notin\{i_1,i_2\}}}_{j<k}}\right.\right.\\
		&\left.\left.+\mathop{\mathop{\sum\sum}_{j,k\notin\{i_1,i_2\}}}_{j<k}(\lambda_{i_1,j,k}+\lambda_{i_2,j,k}+\lambda_{i_1,i_2,j,k})+\cdots+\lambda_{1,2,\cdots,n}\right)\right\}\\
		\vdots&\\
		+&(-1)^{n-1}\exp\left\{-t\left(\sum_{i=1}^n\lambda_i +\mathop{\sum\sum}_{i<j}\lambda_{i,j}+\mathop{\sum\sum\sum}_{i<j<k}\lambda_{i,j,k}+\cdots+\lambda_{1,2,\cdots,n}\right)\right\}\\
	\end{split}
\end{equation*}

$$=\sum_{k=1}^n(-1)^{k-1}\mathop{\sum\cdots\sum}_{1\leq i_1<\cdots<i_k\leq n}\exp\left\{-t\left(\sum_{p=0}^k\left(\mathop{\mathop{\sum\cdots\sum}_{1\leq j_1<\cdots<j_p\leq n}}_{j_1,\cdots,j_p\notin\{i_1,\cdots,i_k\}}\sum_{q=1}^k\left(\mathop{\sum\cdots\sum}_{1\leq u_1<\cdots<u_q\leq k}\lambda_{i_{u_1},\cdots,i_{u_q},j_1,\cdots,j_p}\right)\right)\right)\right\}$$

\subsection{Multivariate Gumbel's Type 1 ($MG1$)}
The resulting joint distribution of lifetimes $X_1,\ldots,X_k$ of the components having $MG1$ is given by
\begin{eqnarray}
\overline{F}(x_1,x_2,\ldots,x_n)&=&\exp\Big\{-\sum_{i=1}^{n}\lambda_i x_i-\sum_{i_1<}\sum_{i_2}\lambda_{i_1,i_2}x_{i_1}x_{i_2}-\sum_{i_1<}\sum_{i_2<}\sum_{i_3}\lambda_{i_1,i_2,i_3}x_{i_1}x_{i_2}x_{i_3}-\ldots-\nonumber\\
&&\lambda_{i_1,i_2,\ldots,i_n}x_{i_1}x_{i_2}\ldots x_{i_n}
\Big\}\nonumber.
\end{eqnarray}
\begin{equation*}
	\begin{split}
		\bar{F}_T(t)=&\sum_{i=1}^n\exp\left\{-t\lambda_i\right\}\\
		&-\mathop{\sum\sum}_{1\leq i_1<i_2\leq n}\exp\left\{-t(\lambda_{i_1}+\lambda_{i_2})-t^2\lambda_{i_1,i_2}\right\}\\
		&+\mathop{\sum\sum\sum}_{1\leq i_1<i_2<i_3\leq n}\exp\left\{-t(\lambda_{i_1}+\lambda_{i_2}+\lambda_{i_3})-t^2(\lambda_{i_1,i_2}+\lambda_{i_1,i_3}+\lambda_{i_2,i_3})-t^3\lambda_{i_1,i_2,i_3}\right\}\\
		&\vdots\\
		&+(-1)^{n-1}\exp\left\{-t\sum_{i=1}^{n}\lambda_i-t^2\mathop{\sum\sum}_{1\leq i_1<i_2\leq n}\lambda_{i_1,i_2}-t^3\mathop{\sum\sum\sum}_{1\leq i_1<i_2<i_3\leq n}\lambda_{i_1,i_2,i_3}-\cdots-t^n\lambda_{i_1,i_2,\ldots,i_n}\right\}\\
		=&\sum_{k=1}^n(-1)^{k-1}\mathop{\sum\cdots\sum}_{1\leq i_1<\cdots<i_k\leq n}\exp\left\{-\sum_{p=1}^k t^p\mathop{\sum\cdots\sum}_{1\leq u_1<\cdots<u_p\leq k}\lambda_{i_{u_1},\cdots,i_{u_p}}\right\}
	\end{split}
\end{equation*}

\end{document}